\documentclass[12pt]{article}
\usepackage{amsthm,amsfonts,amsmath,amssymb,bm}
\usepackage[mathscr]{eucal}
\usepackage[cp1251]{inputenc}
\usepackage[english, russian]{babel}
\usepackage[final]{graphicx}

\newtheorem{theorem}{Теорема}[section]

\begin{document}
	
\begin{center}	
{\bf\Large The full diapason of convergence rates \\ of Birkhoff averages for ergodic flows}

\medskip

{\bf Podvigin I.V., Ryzhikov V.V.}
\end{center}

\text{}

\begin{flushright} 
	\it To Alexander Grigorievich Kachurovskii
	
	on the occasion of his 65th birthday
\end{flushright}

{\bf Abstract:} For an ergodic flow, a range of rates of convergence of Birkhoff averages from the maximum rate to an arbitrarily slow rate is realized by choosing the averaging function. For torus windings, the continuity of the averaging functions is ensured. This complements Krengel's classical result on the slow rates of convergence of means for ergodic automorphisms.

	

\section{Introduction}

For an ergodic measure preserving flow $T_t$ on a probability space $(X,m)$ and a function $f\in L_1(X,m)$, Birkhoff's ergodic theorem states that the time averages
$$
A(f,t,x):=\frac 1 t \int_0^t f(T_s x) ds
$$ 
converge for almost all ${x\in X}$ to the spatial average $\int fdm$.

For some functions, this convergence may be uniform. This type of convergence is the subject of this article. Our goal is to demonstrate how to control the rate of convergence of averages by selecting the function to be averaged. The function is sought as a functional series. As a result, for a given flow, a range of speeds is realized, from the so-called maximum rate to an arbitrarily slow one.	

Continuous realizations of the functions being averaged are found for ergodic torus windings. We also consider the general case of ergodic flows without periodic trajectories. Here, Rudolph's theorem on a special representation of a flow when the return function differs little from a constant is applied. Rudolph's theorem can be viewed as a continuous analogue of the Rokhlin--Halmos lemma. It should be noted that Krengel used this lemma in~\cite{Kr}, and in~\cite{R23}, some delicate generalization of the Rokhlin--Halmos lemma due to Alpern~\cite{A} was used for slow convergences.
The effect of slowing down the convergence rate is also possible in the case of weighted averages, as shown in~\cite{R}. Remarkably, in certain situations, weighted averages yield ultrafast convergence rates of type $o(\frac 1 t)$, see~\cite{Izv}.
Our paper considers only classical Birkhoff averages, for which the maximum rate of convergence of the averages cannot be about $o(\frac 1 t)$ in the case of a nonzero averaged function.	
	
\section{Convergence with maximum rate}
Let $X=[0,1)$. Consider the rotation of a circle as the flow, i.e. $T_t x=\{x+t\}$
(the fractional part of the sum $x+t$). Obviously,
for a function ${f\in L_1(X, m)}$ with zero mean for all $x$, ${A(f,t,x)=0}$ holds for ${t\in\mathbb{Z}}$ and
$$
|A(f,t,x)|\leq \frac {\|f\|_1} {t},\ \ t>0.
$$	
Thus, convergence in Birkhoff's theorem is uniform and occurs at the maximum possible rate (a rate of the form $o(1/t)$ can only occur in the case of a zero function)~\cite{KaPoS2020}. Note that discrete Birkhoff sums
may behave quite differently; see~\cite{Koch23,Koch25} for example.

Above, we mentioned the degenerate case where the entire phase space of the flow is a periodic trajectory. The maximum rate of convergence of averages also arises in the non-ergodic case, when each point of the phase space has a bounded (periodic) trajectory. Specifically, let
$\mathcal{P}(x)$ be the period of the point~${x\in X}$
with respect to the flow ${T_t}$, i.e., the minimum number
$t>0$ such that ${T_tx=x}$; we set for non-periodic points
${\mathcal{P}(x)=\infty}$. If the function ${\mathcal{P}\in L_\infty(X,m)}$, then the convergence of ergodic averages will also be maximum, but, generally speaking, non-uniform. Indeed, assuming $0<t=N\mathcal{P}(x)+r, \, 0\le r<\mathcal{P}(x), N\in\mathbb{N}\cup\{0\},$ we obtain
\begin{multline*}
\left|A(f,t,x)-\frac{1}{\mathcal{P}(x)}\int_0^{\mathcal{P}(x)}f(T_sx)\,ds\right|=
\\
=\left|\frac{1}{t}\sum_{k=0}^{N-1}\int_{k\mathcal{P}(x)}^{(k+1)\mathcal{P}(x)}\!\!f(T_s x)\,ds+
\frac{1}{t}\int_{N\mathcal{P}(x)}^{t}\!\!f(T_s x)\,ds-\frac{1}{\mathcal{P}(x)}\int_0^{\mathcal{P}(x)}\!\!f(T_s x)\,ds\right|=\\
=\left|\left(\frac{N}{t}-\frac{1}{\mathcal{P}(x)}\right)\int_0^{\mathcal{P}(x)}f(T_s x)\,ds+
\frac{1}{t}\int_{0}^{r}f(T_s x)\,ds\right|\leq\\
\leq\frac{2}{t}\int_0^{\mathcal{P}(x)}|f(T_s x)|\,ds.
\end{multline*}	
	
Thus, the presence of slow convergence rates in Birkhoff's ergodic theorem is related to the unboundedness (i.e., non-periodicity) of the flow trajectories. Below, we will consider only such cases.

\vspace{3mm}
{\bf Example of a special ergodic flow.} 
Let $S:[0,1)\to [0,1)$ be a bijection preserving the Lebesgue measure $\mu$ on $[0,1)$. We define a flow $T_t$ on $[0,1)\times [0,1)$ whose action is verbally described as follows: all points $(x,y)$ from ${[0,1)\times [0,1)}$ move vertically upward with unit velocity until they reach the boundary of the square, i.e. ${T_t(x,y)=(x, y+t)}$. Upon reaching the boundary point, ${(x,1)}$ jump to ${(Sx, 0)}$ and then continue moving upward until they end up in $(S^2x,0)$, and so on. We glued
the point ${(x,1)}$ to ${(Sx,0)}$, and the points move with constant velocity,
and this motion preserves area (i.e., the measure $m=\mu\times \mu$ on a square). The resulting flow is ergodic with respect to the measure $m$ only if the transformation $S$ is ergodic with respect to the Lebesgue measure $\mu$ on the interval. Such a flow is called a special flow (or suspension flow) over the automorphism $S$ with constant return function.

Note that for a function independent of $x$ on $[0,1)\times [0,1)$ with zero mean, the assertion of uniform convergence with maximum rate in Birkhoff's ergodic theorem is obvious (the arguments are no different from the case of a rotating circle). Next, we consider a suspension flow over a 2-odometer with a constant return function.

{\bf Dyadic odometer.} Consider a 2-odometer $S$; it operates on the interval $[0,1)$ as follows. For every partition of the interval $[0,1)$ into $2^n$
intervals of the form~${[\frac m n\, ,\frac {m+1} n)}$, the transformation $S$ cyclically permutes these small intervals. Moreover,
each such interval is mapped to itself under the action of the power $S^{2^n}$.
And the restriction of this power to a small interval is similar to the original odometer.

Let us give an explicit definition of such a transformation $S$. We represent the number $x\in [0,1)$ in binary notation: $x=\sum\limits^\infty_{i=1} \frac {x_i} {2^i}$, $x_i\in \{0,1\}$. We set 
$$
S(0,0\dots)=(0,1\dots), \ \ S(0,10\dots)=(0,01\dots),$$ $$S(0,\underbrace{11\dots 11}_{m}0\dots)= (0,\underbrace{00\dots 00}_{m} 1\dots),
$$
where the remaining digits, denoted by ellipses, remain unchanged.
The odometer is also called an adic shift, an addition machine, a von Neumann transform, or an automorphism with a dyadic-rational spectrum. Fig.~1 shows its (self-similar) graph.	

{\bf Special flow over the odometer.} To study the flow over the 2-odometer, for every $n\in\mathbb{N}$ we will identify the space $M=[0,1)\times [0,1)$ with the rectangle ${M_n=[0,\frac 1 {2^n})\times [0, 2^n)}$.

\begin{figure}
	\begin{center}
		\includegraphics[scale=1.2]{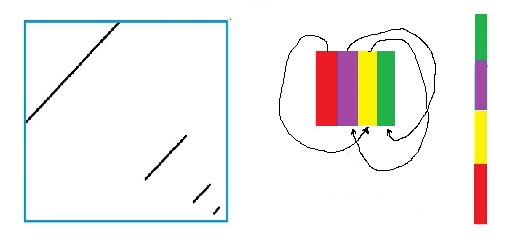} 
		
		{\caption{\small Self-similar graph of 2-odometer (left); Realization of $M_n, n=2$ using the flow over the 2-odometer (right)}}
	\end{center}
\end{figure}

The square $M$ is the union of pieces of trajectories of length~1 for points belonging to the lower boundary. But $M$ can be represented as $M_n$ --- the union of trajectories of length $2^n$ for points from the interval $[0, 2^{-n})\times \{0\}$. In Fig.~1 we depicted this for $n=2$. We introduce coordinates $x_n, y_n$ on such a rectangle $M_n$.

Let us consider functions ${f\in L_1(M)}$ satisfying the following conditions:
\\
(1) ${f(x,y)=\sum\limits_{n=1}^\infty f_n(x_n,y_n)},$ where ${f_n\in L_1(M_n)}$ such that
\\
(2) ${\int_0^{2^n}f_n(x_n,y_n)\, dy_n =0}$ for all $x_n\in[0,2^{-n})$ and 
\\
(3) ${\sum_{n=1}^\infty\|f_n\|_{\infty,1}<\infty},$ where ${\|f_n\|_{\infty,1}=\sup\limits_{x_n\in[0,2^{-n})}\int_0^{2^n}|f_n(x_n, y_n)|\,dy_n}.$

\begin{theorem}\label{Th1}
	Let the function $f$ satisfy conditions (1)-(3). Then the average $A(f,t,(x,y))$ for the flow over the 2-odometer converge uniformly to $0$ as $t\to\infty$ with maximum rate.
\end{theorem}

{\bf Proof.} It is easy to see that condition~(2) implies that the mean value of $f$ is equal to zero. We will show that
for all $t>0$ and all ${(x,y)\in M}$
$$
|A(f,t,(x,y))|\leq\frac{2\sum_{n=1}^\infty\|f_n\|_{\infty,1}}{t}.
$$

Indeed, let ${2^{N-1}<t\leq 2^N},$ then
\begin{multline*}
|A(f,t,(x,y))|\leq\left|A\left(\sum_{n=1}^{N} f_n,\ t,(x,y)\right)\right|+
\left|A\left(\sum_{n=N+1}^{\infty} f_n,\ t,(x,y)\right)\right|\leq\\
\leq\sum_{n=1}^{N}\left|A\left( f_n,\ t,(x,y)\right)\right|+
\sum_{n=N+1}^{\infty }A(|f_n|,\ t,(x,y)).
\end{multline*}

For each term of the second sum we immediately obtain the estimate
\begin{multline*}
A(|f_n|,\ t,(x,y))=\frac{1}{t}\int_0^t|f_n(T_s(x_n,y_n))|\,ds=\\
=\frac{1}{t}\int_{0}^{2^n-y_n}|f_n(x_n,s+y_n)|\,ds+\frac{1}{t}\int_0^{t+y_n-2^n}|f_n(x'_n,s)|\,ds\leq\\
\leq\frac{1}{t}\int_{0}^{2^n}|f_n(x_n,s)|\,ds+\frac{1}{t}\int_0^{2^n}|f_n(x'_n,s)|\,ds\leq
\frac{2\|f_n\|_{\infty,1}}{t}.
\end{multline*}
Here we assumed that ${t+y_n>2^n},$ otherwise the factor 2 would not appear in the final estimate.

For the terms of the first sum we use~condition (2):
\begin{multline*}
tA(f_n,\ t,(x,y))=\int_0^{2^n-y_n}f_n(x_n,y_n+s)\,ds+
\underbrace{\int_0^{2^n}f_n(x^{(1)}_n,s)\,ds}_{=0}+...\\
...+\underbrace{\int_0^{2^n}f_n(x^{(k)}_n,s)\,ds}_{=0}+\int_0^{t-k2^n-y_n}f_n(x^{(k+1)}_n,s)\,ds.
\end{multline*}
Here $k\in\mathbb{N}$ is the largest natural number for which ${t-2^nk-y_n\geq0}.$ From this we obtain an estimate for the terms of the first sum
$$
|A(f_n,t, (x,y))|\leq\frac{1}{t}\left|\int_0^{2^n-y_n}f_n(x_n,y_n+s)\,ds+\int_0^{t-k2^n-y_n}f_n(x^{k+1}_n,s)\,ds\right|\leq
\frac{2\|f_n\|_{\infty,1}}{t}.
$$
Summing the inequalities for both sums, we obtain the required estimate for the time averages. Theorem~\ref{Th1} is proved.

\section{Slow convergence of Birkhoff averages for flow over a dyadic odometer}

Removing condition~(3) from the above example allows us to obtain an arbitrarily slow rate of convergence of the time averages. Let us demonstrate this.

Let $p_n$ be a rapidly growing sequence of natural numbers. Denote $L_n=d_n2^{p_n-1}$, where ${d_n\in(0,1/2)}$. Consider on ${M_n=[0, 2^{-p_n})\times [0, 2^{p_n})}$ a function $f_n$, which on the rectangle $[0, 2^{-p_n})\times [0, L_n)$ is equal to ${a_n>0}$, on the rectangle ${[0, 2^{-p_n})\times [2^{p_n-1}, 2^{p_n-1}+L_n)}$ is equal to $-a_n$, and at the remaining points $M_n$ is equal to zero. It is easy to verify that for any ${x_n\in[0,2^{-p_n})}$
$$
\int_0^{2^{p_n}}f_n(x_n,y_n)\,dy_n=0,\ \ \int_0^{2^{p_n}}|f_n(x_n,y_n)|\,dy_n=2a_nL_n.
$$
To ensure that condition~(3) of Lemma~1 is not satisfied for the function ${f(x,y)=\sum\limits_{n=1}^\infty f_n(x_n,y_n)}$, we impose the condition
$$
\sum_{n=1}^{\infty}a_nL_n=\sum_{n=1}^{\infty}a_nd_n2^{p_n}=\infty.
$$
Moreover, for $f$ to belong to ${L_1(M)}$, it is sufficient that the series of $L_1(M_n)$-norms of the functions $f_n$ converge, i.e. series ${\sum\limits_{n=1}^{\infty}a_nd_n}.$ We will require a stronger condition
$$
\sum_{m=n+1}^\infty a_m=o(a_nd_n)\ \ \text{for}\ \ n\to+\infty.
$$

\begin{theorem}\label{Th2}
	Let $f$ be the function constructed above and $t_n=2^{p_n-2}.$ Then the sequence of averages $A(f,t_n,(x,y))$ for the flow over the 2-odometer converges uniformly to $0$ as $n\to\infty$ with a rate of $\mathcal{O}(a_nd_n)$.
\end{theorem}

{\bf Proof.}
Representing
$$
A(f, t_n, (x,y))=A\left(\sum_{m=1}^{n-1}f_m, t_n, (x,y)\right)+A(f_n, t_n, (x,y))+
A\left(\sum_{m=n+1}^{\infty}f_m, t_n, (x,y)\right),
$$
we note that for all ${(x,y)\in M}$
$$
A\left(\sum_{m=1}^{n-1}\, f_m, t_n, (x,y)\right)=0.
$$
This follows from the fact that the integration interval $[0,t_n)$ is divided into a finite number of intervals of length $2^{p_m}$ for each ${m=1,...,n-1}.$ This will indeed be the case, since $2^{p_n-2}$ is divisible by $2^{p_m}$. This leads to a more precise condition on the increasing sequence $p_n:$
$$
p_{n+1}\geq p_n+2.
$$
And on an interval of length $2^{p_m}$, the integral of $f_m$ with respect to the variable $y_m$ is equal to zero.

Now let's look at the distribution of the values of the function
$$
{A_n(x,y)=A\left(f_n, t_n, (x,y)\right)}.
$$
On a set of measure ${\frac 1 2 - d_n}$, the function $A_n(x,y)$ is equal to $0$, on a set of measure ${1/ 4 - d_n/2}$, it is equal to $2a_nd_n$ (the maximum value), and on a set of the same measure, $A_n(x,y)$ is equal to $-2a_nd_n$ (the minimum value).
The remaining intermediate values vary linearly.

Let us now show how the values of the function $A\left(\sum_{m=1}^\infty \, f_m, t_n, (x,y)\right)$ are estimated. We will assume that for all $m>n$
the condition
$$
2^{p_n-2}<d_m2^{p_m-1}
$$ 
is satisfied. Since the length of the integration interval, equal to $2^{p_n-2},$
is now significantly smaller than $L_m$ for $m>n,$ the largest absolute value of the mean $A(f_m, t_n, (x,y))$ will be $a_m.$ Thus,
$$
\left|A\left(\sum_{m=n+1}^\infty\, f_m, t_n , (x,y)\right)\right|\leq\sum_{m=n+1}^\infty\,A(|f_m|, t_n , (x,y))\leq\sum_{m=n+1}^\infty a_m=o(a_nd_n).
$$
Putting the estimates together, we obtain for all ${(x,y)\in M}$
$$
A(f, t_n, (x,y))=\mathcal{O}(a_nd_n)\ \ \text{as}\ \ n\to\infty.
$$
Theorem~\ref{Th2} is proved.

The values of $A\left(f, t_n ,(x,y)\right)$ on a set of measure close to $1/2$ are close in absolute value to $2a_nd_n$. And for most of the remaining points, the values of $A\left(f, t_n ,(x,y)\right)$ are asymptotically $o(a_nd_n)$. Thus, it is impossible to obtain an estimate of the form $o(a_nd_n).$ Moreover, for given sequences $a_n$ and $d_n$, we can choose $t_n\to\infty$ growing arbitrarily quickly.
Thus, we obtain arbitrarily slow convergence of the Birkhoff averages, and the distribution of the values of these means is such that over almost half the space, the function takes values asymptotically infinitely large compared to the values it takes over the remaining part of the space.

\section{Slow convergence of Birkhoff averages in the general case}

We have shown how slow convergence of means for a flow over a 2-odometer is realized. A similar effect can be obtained in the general case. We partially repeat the previous construction, using the Rudolph special rectangular representation of the flow and Birkhoff's ergodic theorem when choosing $t_{n+1}$ and $f_{n+1}$. The moment $t_{n+1}$ is chosen such that the means $A\left(f_1+\dots +f_n, t_{n+1} ,x\right)$ for most $x$ are extremely small compared to $a_{n+1}$.

To choose $f_{n+1}$, we use the theorem on a special rectangular representation~(see, for example,~\cite[Chapter~11,\S4]{KSF}). By Rudolph's theorem, for an aperiodic ergodic flow $T_t$ on the probability space~${(X,m)}$ and any positive numbers $p,q$ such that $p/q$ is irrational, there exists a special representation of the flow with a function taking two values: $p$ and $q$. Thus, the phase space of the flow can be partitioned into two measurable sets, which we identify with rectangles. The first set consists of trajectory segments of length $p$, and the second set consists of segments of length $q$.

\begin{theorem}\label{Th4}
	Let $T_t$ be an aperiodic ergodic flow on the probability space~${(X,m)}$ and ${\varphi(t)\to+0}.$ Then there exists a function ${f\in L_1(X,m)}$ with zero mean such that for almost all ${x\in X}$
	$$
	\limsup_{t\to+\infty}\frac{1}{\varphi(t)}|A(f,t,x)|=+\infty.
	$$
\end{theorem}

{\bf Proof.} Without loss of generality, we can assume that $\varphi(t)$ monotonically tends to zero. We seek the function ${f\in L_1(X,m)}$ in the form
$$
f(x)=\sum_{n=1}^\infty f_n(x_n,y_n),
$$
where $f_n$ is constructed as follows. In the Rudolph representation,
we choose the height parameters of the rectangles: $p=h_n$ and $q=h_n+\varepsilon_n$. Moreover,
$$
h_nc_n+(h_n+\varepsilon_n)d_n=1,
$$
and ${h_n\to+\infty},$ ${\varepsilon_n\to0}.$

The function $f_n(x_n,y_n)$ will depend only on the height of $y_n$.
Specifically, let $f_n=a_n$ for $0\leq y_n<h_n/4$; $f_n=-a_n$ for $h_n/2\leq y_n< 3h_n/4$, and let $f_n$ be 0 elsewhere (see Figure~2).

\begin{figure}
	\begin{center}
		\includegraphics[scale=0.5]{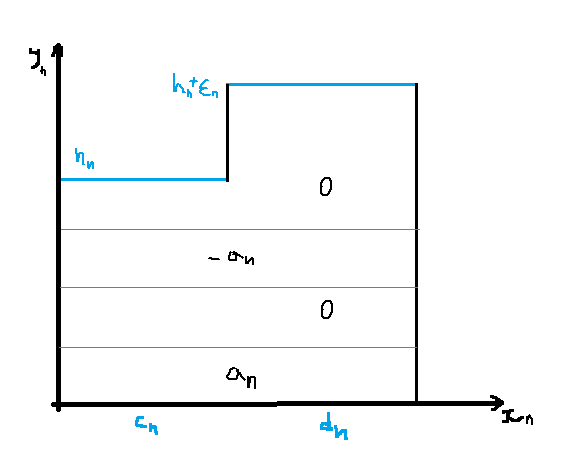}
		
		{\caption{\small Rudolph representation with parameters $h_n$ and $h_n+\varepsilon_n$}}
	\end{center}
\end{figure}

We can measurably realize this entire "picture"\ in the original space ${(X,m)}$ by an inverse isomorphism (the point $x$ corresponds to the pair $(x_n,y_n)$). We leave the notation for the function the same. Let ${\delta_n\to+0}$ and ${t_n=h_n\delta_n}.$ To control the ergodic averages, we define asymptotic conditions on all sequences:
$$
\sum_{n=1}^\infty\varepsilon_nd_n+4\delta_n<1;\eqno(I)
$$
$$
a:=\sum_{n=1}^\infty a_n<+\infty\ \ \text{and}\ \ \sum_{m=n+1}^\infty a_m=\alpha_na_n=o(a_n);\eqno(II)
$$  
$$
\frac{h_1+h_2+...+h_{n-1}}{h_na_n\delta_n}\to0\ \ \ \text{and}\ \ \ \frac{a_n\alpha_n}{\varphi(t_n)}\to+\infty.\eqno(III)
$$

By condition $(I)$, we define the sequences ${\varepsilon_n}$ and ${\delta_n}$. Next, by condition~$(II)$, we define $a_n$ and, most importantly for controlling ergodic averages, we define the growth $h_n$ by conditions $(III).$

As before, we represent $A(f,t_n,x)$ as the sum of three terms:
$$
A(f, t_n, x)=A\left(\sum_{m=1}^{n-1}f_m, t_n, x\right)+A(f_n, t_n, x)+
A\left(\sum_{m=n+1}^{\infty}f_m, t_n, x\right).
$$
For the third term from condition~$(II)$ for all ${x\in X}$ we have
$$
\left|A\left(\sum_{m=n+1}^{\infty}f_m, t_n, x\right)\right|\leq\sum_{m=n+1}^\infty a_m=o(a_n).
$$
For the first term from the first part of condition~$(III)$, in view of the vanishing of the integral on entire segments of the trajectories, for all ${x\in X}$ we have
$$
A\left(\sum_{m=1}^{n-1}|f_m|, t_n, x\right)\leq\frac{2}{t_n}\sum_{m=1}^{n-1}a_mh_m\leq
2a\frac{h_1+h_2+...+h_{n-1}}{h_n\delta_n}=o(a_n).
$$
For the second term, it's easy to calculate that it takes the value $a_n$ on a set of measure $(c_n+d_n)h_n(1/4-\delta_n)$, the value $-a_n$ on another set of the same measure, and zero on a set of measure at least $2(c_n+d_n)h_n(1/4-\delta_n)$. The remaining values do not exceed $a_n$ in absolute value.

Thus, $A(f,t_n, x)=\mathcal{O}(a_n)$ for all $x\in X.$ Moreover, on some set $\mathcal D_n$ of measure
$$
m(\mathcal{D}_n)=2(c_n+d_n)h_n(1/4-\delta_n)=2(1/4-\delta_n)(1-\varepsilon_nd_n)>1/2-2\delta_n-\frac{\varepsilon_nd_n}{2}
$$
$$
|A(f,t_n,x)|=a_n+o(a_n);
$$
and on some set $\mathcal{E}_n,$ with the same estimate for the measure
$$
m(\mathcal{E}_n)>1/2-2\delta_n-\frac{\varepsilon_nd_n}{2},
$$
there will be
$$
|A(f,t_n,x)|=o(a_n).
$$
Let ${\mathcal{F}=\bigcap_{n=1}^\infty(\mathcal{E}_n\cup\mathcal{D}_n)},$ then, given~$(I),$ ${m(\mathcal{F})>1-\sum_{n=1}^\infty\varepsilon_nd_n+4\delta_n>0}.$

For all $x\in\mathcal{F}$, for the averages under consideration, we have
$$
\text{either}\ \ |A(f,t_n,x)|=a_n(1+\alpha_n), \ \ \text{or}\ \ |A(f,t_n,x)|=\alpha_na_n.
$$
In any case, for such $x$, from the second condition in~$(III)$ we obtain
$$
\frac{ |A(f,t_n,x)|}{\varphi(t_n)}\to+\infty,
$$
i.e., for $x\in\mathcal{F}$
$$
\limsup_{t\to+\infty}\frac{1}{\varphi(t)}|A(f,t,x)|=+\infty.
$$
For a.e. $x\in X\setminus \mathcal{F}$ the same relation will hold due to the zero-one law for the rate of convergence~\cite{KaPoS2022}.

\section {On the slow convergence of time averages for a torus winding and a continuous function}

It turns out that averaging a continuous function along the trajectories of a smooth flow can be combined with the effect of slow convergence of averages. To do this, we will smooth the functions $f_n$ considered above in a suitable manner. Let $T_t(x,y)=(\{x+t\},\{y+ct\}), {(x,y)\in M}$ be an ergodic winding of a torus along
the vector $(1,c)$, where $c$ is an irrational number approximated by rational fractions $p_n/q_n$ such that
$$
0<c-\frac{p_n}{q_n}=\beta_n\frac{1}{q_n^2}=o\left(\frac{1}{q_n^2}\right).
$$

Considering the winding of a torus along the vector ${(1,p_n/q_n)}$, we identify the torus~$M$ with a narrow parallelogram~$M_n$ with a side of length $\sqrt{p_n^2+q_n^2}$ parallel to this vector, and a side of length $1/q_n$ parallel to the OX axis. Inside this parallelogram, consider a parallelogram~$R_n$ directed along the flow $T_t.$
\begin{figure}
	\begin{center}
		\includegraphics[scale=0.5]{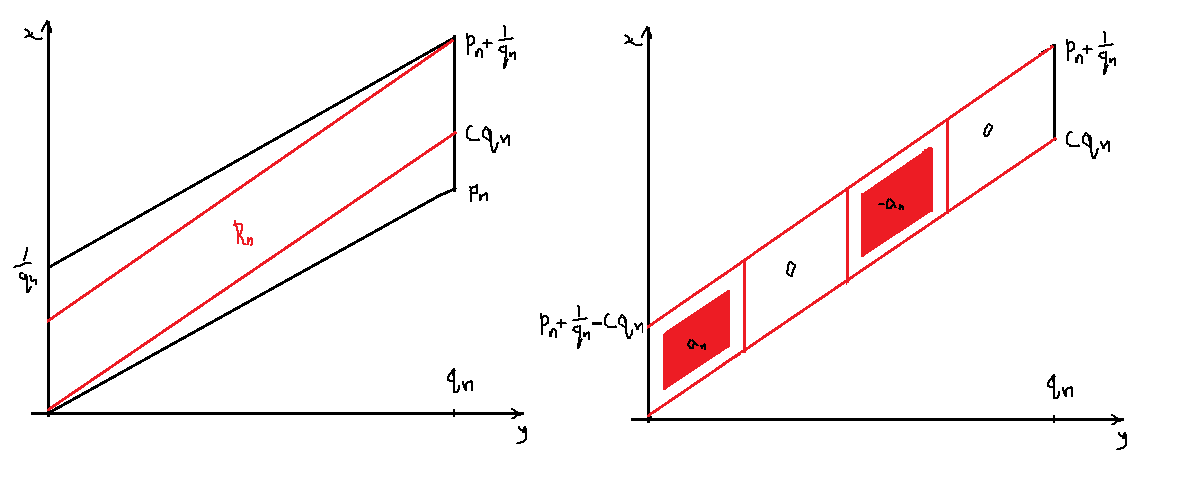}
		
		{\caption{\small Parallelogram $R_n$ and definition of function $f_n$ on it}}
	\end{center}
\end{figure}

It is easy to see that the area of the parallelogram $R_n$ is
$$
1-q_n^2\left(c-\frac{p_n}{q_n}\right)=1-\beta_n.
$$

On $M_n$, we choose coordinates $x_n\in[0,1/q_n]$ along the OX axis and $y_n\in[0,q_n\sqrt{1+c^2}]$ along the vector $(1,c).$
We construct a smooth function $f_n(x_n,y_n)$ whose absolute value does not exceed ${a_n>0}$, such that ${\sum_n a_n<\infty}$. Specifically, we divide $R_n$ into four parts of equal area. In the first quadrant, $f_n$ is equal to $a_n$ on the set separated from the boundary and an area close to $1/4$. In the remaining part of this quadrant, $f_n$ decreases smoothly to zero at the boundary. In the second and fourth parallelograms, $f_n$ is equal to 0. In the third quadrant, the values are antisymmetric to the values in the first quadrant (see Fig.~3).

Outside the rectangle $R_n$ and on its boundary, let $f_n$ equal $0$. We also assume that
$$
\int_0^{q_n\sqrt{1+c^2}}f_n(x_n,y_n)\,dy_n=0
$$
for each ${x_n\in[0,1/q_n]}.$
As a result, we obtain the continuous function $f(x,y)=\sum\limits_{n=1}^\infty f_n(x_n,y_n)$.

Let ${\varphi(t)\to+0}$ as $t\to+\infty.$ Set $t_n=\delta_nq_n\sqrt{1+c^2},$ and let conditions $(II)$ and $(III)$ be satisfied, where $h_n$ is replaced by $q_n.$

\begin{theorem}\label{Th3}
	Let $f$ be the function constructed above. Then the ergodic averages $A(f,t,(x,y))$ for the torus winding along the vector ${(1,c)}$ satisfy for almost all ${(x,y)\in M}$ the relation
	$$
	\limsup_{t\to+\infty}\frac{1}{\varphi(t)}|A(f,t,(x,y))|=+\infty.
	$$
\end{theorem}

{\bf Proof.} We represent the time average as
$$
A(f,t_n,(x,y))=A\left(\sum_{m=1}^{n-1}f_m, t_n, (x,y)\right)+A(f_n,t_n,(x,y))+A\left(\sum_{m=n+1}^\infty f_m, t_n, (x,y)\right).
$$
For all ${(x,y)\in M}$, we have an obvious estimate for the third term:
$$
\left|A\left(\sum_{m=n+1}^\infty f_m, t_n, (x,y)\right)\right|\leq \sum_{m=n+1}^nA(|f_m|, t_n, (x,y))\leq\sum_{m=n+1}^\infty a_m=o(a_n).
$$
To estimate the first term in the above sum, we take advantage of the vanishing integrals along the flow trajectory. On most of the trajectory segments along which the integral is taken, the contribution to the integral is zero due to the choice of the functions $f_m$.
We only need to estimate the values of the integrals along the initial and final parts of the trajectory segments. Namely, for each $m=1,...,n$ and all $(x,y)\in M$
\begin{multline*}
t_nA(f_m, t_n, (x,y))=\int_0^{q_m\sqrt{1+c^2}-y_m}f_m(T_s(x_m,y_m))\,ds+
\underbrace{\int_0^{q_m\sqrt{1+c^2}}f_m(x^{(1)}_m,s)\,ds}_{=0}+...\\
...+\underbrace{\int_0^{q_m\sqrt{1+c^2}}f_m(x^{(k-1)}_m,s)\,ds}_{=0}+\int_0^{y_m}f_m(x^{(k)}_m,s)\,ds.
\end{multline*}
It follows
\begin{multline*}
\left|A\left(\sum_{m=1}^{n-1}f_m, t_n, (x,y)\right)\right|\leq\sum_{m=1}^{n-1}\frac{2a_mq_m\sqrt{1+c^2}}{t_n}=\mathcal{O}\left(\frac{q_1+\dots+q_{n-1}}{\delta_nq_n}\right)=o(a_n).
\end{multline*}

The second term behaves the same as in Theorem~\ref{Th4} for the general case. Namely, on a set close in measure to $1/2$, the averages ${|A(f,t_n,(x,y))|=a_n+o(a_n)},$ and on another set, also close in measure to $1/2$, the averages behave like ${o(a_n)}.$
Theorem~\ref{Th3} is proved.

\newpage
\section{Concluding remarks}

The topic of convergence rates of ergodic averages is associated with a large number of unsolved problems.

{\bf Convergence of averages and smoothness of the averaged function.} The function constructed in Theorem~\ref{Th3} is only continuous. Of interest is the question of the existence of smooth functions with an arbitrarily slow rate of convergence of their averages. As Kowada~\cite{Kow} (see also~\cite{Web}) showed, starting from a certain smoothness index depending on the rate of approximation of the irrational number $c$, the rate of convergence of ergodic averages for a torus winding will be maximal, i.e., $\mathcal{O}(1/t).$

\vspace{3mm}

{\bf Possible distributions of Birkhoff averages.}
We have shown how the choice of an appropriate averaging function realizes a range of convergence rates from the maximum to arbitrarily slow. Slowing effects on the rate of convergence of averages can be detected for a wide class of group actions~(see~\cite{R}). Of interest is not only the estimation of convergence rates, but also the more general problem of possible distributions of Birkhoff average values. By choosing the averaging function, these distributions can be controlled. For example, suppose we want the distribution of the values of the function
$A(f, t_n, x)$ for an ergodic flow to be arbitrarily close to the distribution of the values of, say, the function $c_n s^2 , s\in[0, 1],$ with respect to Lebesgue measure on $[0, 1].$ If the numbers $c_n$ tend to $+0$ sufficiently quickly, the problem can be solved by a method similar to the one described above, with an appropriate choice of the functions $f_m$.

\vspace{3mm}

{\bf Optimal weight distributions.} As already noted, when considering weighted averages, it is possible to increase the range of rapid rates of convergence for averages (see~\cite{Izv}, as well as references therein). This raises a number of new problems. We formulate the following special cases.

Let $f$ be a zero-mean function of given smoothness $k\geq1$. For an ergodic shift $T$ on the torus, consider all possible convex sums of the form
$$
A_w(f,N,x)=\sum _{n=1}^Nw_{n,N}f(T^nx),\ \
$$
where $w=(w_{1,N},\dots,w_{N,N})$ is a probability vector.

For what probability distributions of the coefficients $w_{n,N}$, given a sufficiently large $N$, will the norms $\|A_w\|$ (in $L_1$, $L_2$, or $L_\infty$ spaces) be minimal?

A similar question is posed for continuous time and torus windings.
For average
$$
A_w(f,t,x) =\int_0^t w(s,t) f(T_s x)\,ds
$$
where $w\geq 0$ and $\int_0^t w(s,t)ds =1$ for each $t>0,$
we seek $\inf_{w}\|A_w\|.$

Even in the case of a rotating circle, such problems appear nontrivial.

\newpage

\begin{center}	
{\bf\Large Полный диапазон скоростей сходимости \\ средних Биркгофа  для эргодических  потоков}

\medskip

{\bf  Подвигин И.В., Рыжиков В.В.}
\end{center}

\begin{flushright} 
\it Качуровскому Александру Григорьевичу 
	
по случаю его 65-летия
\end{flushright}

{\bf Аннотация:} 

Для  эргодического потока реализуется 
диапазон скоростей сходимости средних Биркгофа от максимальной скорости до сколь угодно медленной путем подбора усредняемой функции.  Для обмоток тора при этом обеспечена   непрерывность  усредняемых функций.    
Это дополняет  классический  результат Кренгеля о медленных скоростях сходимости средних для эргодических автоморфизмов.  

{\bf Ключевые слова:} скорости сходимости в эргодических теоремах, специальное представление потока, одометр, обмотка тора

\setcounter{section}0
\setcounter{figure}0

\section{Введение}

Для эргодического сохраняющего меру потока $T_t$ на вероятностном пространстве $(X,m)$ и функции $f\in L_1(X,m)$ теорема Биркгофа утверждает, что  временн\'{ы}е средние 
$$
A(f,t,x):=\frac 1 t \int_0^t f(T_s x) ds
$$ сходятся для почти всех ${x\in X}$  к пространственному среднему $\int fdm$. 

Для некоторых функций указанная сходимость может быть  равномерной. Именно такая сходимость  фигурирует  в настоящей статье.
Наша цель  --- показать управление скоростью сходимости средних
путем подбора  усредняемой  функции. Последняя ищется в виде  функционального ряда. В результате для заданного потока реализуется 
диапазон скоростей от так называемой максимальной скорости до сколь угодно медленной.  
Для эргодических обмоток тора при этом найдены непрерывные реализации усредняемых функций. Мы также рассмотрим  общий случай эргодических потоков без периодических траекторий.  Здесь применяется теорема  Рудольфа о специальном представлении потока, когда функция возвращения  мало отличается от константы. Теорему Рудольфа  можно рассматривать как  непрерывный аналог  леммы Рохлина--Халмоша. Следует отметить, что Кренгель в работе~\cite{Kr}  использовал эту лемму,
а в~\cite{R23} для медленных сходимостей применялось деликатное обобщение леммы Рохлина--Халмоша, принадлежащее Альперну~\cite{A}.
Эффект замедления скорости сходимости возможен и в случае весовых усреднений, что показано в~\cite{R}. Замечательно, что в определенных ситуациях весовые усреднения дают сверхбыстрые сходимости со скоростью~$o(\frac 1 t)$, см.~\cite{Izv}.   
В предлагаемой статье рассматриваются только классические  средние Биркгофа, для которых  максимальная скорость сходимости средних не может быть  $o(\frac 1 t)$ в случае  ненулевой усредняемой функции. 

\section{Сходимоcть с максимальной скоростью}
Пусть $X=[0,1)$, в качестве потока рассмотрим вращение окружности, т.е.  $T_t x=\{x+t\}$
(дробная часть суммы $x+t$).  Очевидно,
что для функции ${f\in L_1(X, m)}$ с нулевым средним для всех $x$  выполнено ${A(f,t,x)=0}$ при ${t\in\mathbb{Z}}$ и 
$$
|A(f,t,x)|\leq  \frac {\|f\|_1} {t},\ \  t>0.
$$
Таким образом, сходимость в теореме Биркгофа равномерная с максимальной возможной скоростью (скорость вида $o(1/t)$ может быть только в случае нулевой функции)~\cite{KaPoS2020}. Отметим, что дискретные суммы Биркгофа
могут вести себя совсем иначе, см., например,~\cite{Koch23,Koch25}.

Выше мы упомянули вырожденный случай,  когда все фазовое пространство потока  является  периодической траекторией. Максимальная скорость сходимости средних  возникает и в неэргодическом случае, когда каждая точка фазового пространства имеет ограниченную (периодическую) траекторию. А именно, пусть
$\mathcal{P}(x)$ --- период точки~${x\in X}$
относительно потока ${T_t},$ т.е. такое минимальное число
$t>0,$ что ${T_tx=x};$ полагаем для непериодических точек
${\mathcal{P}(x)=\infty}$. Если функция ${\mathcal{P}\in L_\infty(X,m)},$ то сходимость эргодических средних также будет с максимальной скоростью, но уже, вообще говоря, не равномерная. Действительно, полагая $0<t=N\mathcal{P}(x)+r, \, 0\le r<\mathcal{P}(x), N\in\mathbb{N}\cup\{0\},$ получаем
\begin{multline*}
\left|A(f,t,x)-\frac{1}{\mathcal{P}(x)}\int_0^{\mathcal{P}(x)}f(T_sx)\,ds\right|=
\\
=\left|\frac{1}{t}\sum_{k=0}^{N-1}\int_{k\mathcal{P}(x)}^{(k+1)\mathcal{P}(x)}\!\!f(T_s x)\,ds+
\frac{1}{t}\int_{N\mathcal{P}(x)}^{t}\!\!f(T_s x)\,ds-\frac{1}{\mathcal{P}(x)}\int_0^{\mathcal{P}(x)}\!\!f(T_s x)\,ds\right|=\\
=\left|\left(\frac{N}{t}-\frac{1}{\mathcal{P}(x)}\right)\int_0^{\mathcal{P}(x)}f(T_s x)\,ds+
\frac{1}{t}\int_{0}^{r}f(T_s x)\,ds\right|\leq\\
\leq\frac{2}{t}\int_0^{\mathcal{P}(x)}|f(T_s x)|\,ds.
\end{multline*}

Таким образом, наличие  медленных скоростей сходимости в эргодической теореме Биркгофа связано с неограниченностью (т.е. непериодичностью) траекторий потока. Ниже мы будем рассматривать только такие случаи. 

\vspace{3mm} 
\bf Пример специального  эргодического потока. \rm   
Пусть $S:[0,1)\to [0,1)$ --- биекция, сохраняющая меру Лебега $\mu$ на $[0,1)$. Определим  поток $T_t$ на $[0,1)\times [0,1)$, действие которого вербально описывается так: все точки $(x,y)$  из ${[0,1)\times [0,1)}$ движутся с единичной скоростью вертикально вверх до встречи с границей квадрата, т.е. ${T_t(x,y)=(x, y+t)}$.  При достижении граничной точки ${(x,1)}$ перескакивают в ${(Sx, 0)}$ и далее продолжают движение вверх, пока не окажутся в $(S^2x,0)$ и т.д. Мы склеили 
точку ${(x,1)}$ с ${(Sx,0)},$ и точки движутся с постоянной скоростью, 
и при этом движении сохраняется площадь (т.е. мера $m=\mu\times \mu$ на квадрате). Полученный поток эргодичен относительно меры $m$ только в случае эргодичности преобразования $S$ относительно меры Лебега $\mu$ на отрезке. Такой поток  называют специальным потоком над автоморфизмом $S$ с постоянной функцией возвращения.

Заметим, что для независящей от $x$  функции на $[0,1)\times [0,1)$ с нулевым средним утверждение о равномерной сходимости с максимальной скоростью в теореме Биркгофа  очевидно (аргументация не отличаются от случая вращения окружности).  Далее рассматривается  специальный поток над 2-одометром с постоянной функцией возвращения.  

\bf  Диадический одометр. \rm Рассмотрим 2-одометр $S$, он действует на отрезке  $[0,1)$ следующим образом. Для всякого разбиения отрезка $[0,1)$ на $2^n$
отрезков  вида~${[\frac m n\, ,\frac {m+1}  n)}$  преобразование $S$ циклически переставляет эти маленькие отрезки. При этом  
каждый такой отрезок переходит в себя при действии степени $S^{2^n}$.
А ограничение этой степени на маленький отрезок подобно исходному одометру.

Дадим явное определение такого  преобразования $S$.  Число $x\in [0,1)$ представим в двоичной записи:  $x=\sum\limits^\infty_{i=1} \frac {x_i} {2^i}$, $x_i\in \{0,1\}$.  Положим  $$S(0,0\dots)=(0,1\dots), \ \ S(0,10\dots)=(0,01\dots),$$   $$S(0,\underbrace{11\dots 11}_{m}0\dots)= (0,\underbrace{00\dots 00}_{m} 1\dots),$$
где остальные разряды, обозначенные многоточием, остаются без изменений.  
Одометр также называется адическим сдвигом, машиной сложения, преобразованием фон Неймана, автоморфизмом с двоично-рациональным спектром. На Рис.~1 изображен его  (самоподобный) график.

\bf Специальный поток над одометром. \rm  Для изучения  потока над 2-одометром мы будем для всякого $n\in\mathbb{N}$  отождествлять пространство $M=[0,1)\times [0,1)$ с прямоугольником ${M_n=[0,\frac 1 {2^n})\times [0, 2^n)}$.
\begin{figure}
\begin{center}
	\includegraphics[scale=1.2]{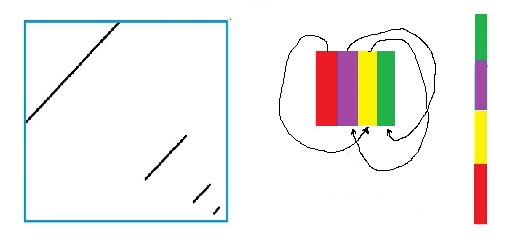} 
	
	{\caption{\small Самоподобный график 2-одометра (слева); реализация $M_n, n=2$ с помощью потока над 2-одометром (справа)}}
\end{center}
\end{figure}
Квадрат $M$ есть  объединение кусков траекторий длины~1 для точек, принадлежащих нижней границе. Но $M$ можно представить как $M_n$ --- объединение траекторий длины $2^n$ для точек из отрезка $[0, 2^{-n})\times \{0\}$. На Рис.~1 мы изобразили это для $n=2$.   Введем координаты $x_n, y_n$ на таком прямоугольнике $M_n$. 

Рассмотрим функции ${f\in L_1(M)}$ удовлетворяющие следующим условиям:
\\
(1) ${f(x,y)=\sum\limits_{n=1}^\infty f_n(x_n,y_n)},$ где ${f_n\in L_1(M_n)}$ такие, что 
\\
(2) ${\int_0^{2^n}f_n(x_n,y_n)\, dy_n =0}$ для каждого $x_n\in[0,2^{-n})$ и 
\\
(3) ${\sum_{n=1}^\infty\|f_n\|_{\infty,1}<\infty},$ где ${\|f_n\|_{\infty,1}=\sup\limits_{x_n\in[0,2^{-n})}\int_0^{2^n}|f_n(x_n, y_n)|\,dy_n}.$

\begin{theorem}\label{Th1} 
Пусть функция $f$ удовлетворяет условиям (1)-(3). Тогда средние $A(f,t,(x,y))$ для потока над 2-одометром равномерно сходятся к $0$ при $t\to\infty$ с максимальной скоростью.
\end{theorem}

{\bf Доказательство.} Нетрудно видеть, что условие~(2) влечет за собой равенство нулю среднего значения функции $f.$ Покажем, что
для всех $t>0$ и всех ${(x,y)\in M}$
$$
|A(f,t,(x,y))|\leq\frac{2\sum_{n=1}^\infty\|f_n\|_{\infty,1}}{t}.
$$
		
Действительно, пусть ${2^{N-1}<t\leq 2^N},$ тогда

\begin{multline*}
|A(f,t,(x,y))|\leq\left|A\left(\sum_{n=1}^{N} f_n,\ t,(x,y)\right)\right|+
\left|A\left(\sum_{n=N+1}^{\infty} f_n,\ t,(x,y)\right)\right|\leq\\
\leq\sum_{n=1}^{N}\left|A\left( f_n,\ t,(x,y)\right)\right|+
\sum_{n=N+1}^{\infty }A(|f_n|,\ t,(x,y)).
\end{multline*}

Для каждого слагаемого второй суммы сразу получаем оценку
\begin{multline*}
A(|f_n|,\ t,(x,y))=\frac{1}{t}\int_0^t|f_n(T_s(x_n,y_n))|\,ds=\\
=\frac{1}{t}\int_{0}^{2^n-y_n}|f_n(x_n,s+y_n)|\,ds+\frac{1}{t}\int_0^{t+y_n-2^n}|f_n(x'_n,s)|\,ds\leq\\
\leq\frac{1}{t}\int_{0}^{2^n}|f_n(x_n,s)|\,ds+\frac{1}{t}\int_0^{2^n}|f_n(x'_n,s)|\,ds\leq
\frac{2\|f_n\|_{\infty,1}}{t}.
\end{multline*}
Здесь мы считали, что ${t+y_n>2^n},$ в противном случае не возникнет множителя 2 в итоговой оценке.

Для слагаемых первой суммы воспользуемся условием (2):
\begin{multline*}
tA(f_n,\ t,(x,y))=\int_0^{2^n-y_n}f_n(x_n,y_n+s)\,ds+
\underbrace{\int_0^{2^n}f_n(x^{(1)}_n,s)\,ds}_{=0}+...\\
...+\underbrace{\int_0^{2^n}f_n(x^{(k)}_n,s)\,ds}_{=0}+\int_0^{t-k2^n-y_n}f_n(x^{(k+1)}_n,s)\,ds.
\end{multline*}
Здесь $k\in\mathbb{N}$ есть наибольшее натуральное число, для которого ${t-2^nk-y_n\geq0}.$ Отсюда уже получаем оценку для слагаемых первой суммы
$$
|A(f_n,t, (x,y))|\leq\frac{1}{t}\left|\int_0^{2^n-y_n}f_n(x_n,y_n+s)\,ds+\int_0^{t-k2^n-y_n}f_n(x^{k+1}_n,s)\,ds\right|\leq
\frac{2\|f_n\|_{\infty,1}}{t}.
$$
Суммируя неравенства для обеих сумм, получим требуемую оценку для временных средних. Теорема~\ref{Th1} доказана.

\section{Медленная сходимость средних Биркгофа для потока над  диадическим одометром}

Отказ от условия~(3) из выше изложенного  примера  позволяет  получить сколь угодно медленную скорость сходимости временных средних. Покажем это.

Пусть задана быстро растущая последовательность $p_n$ натуральных чисел. Обозначим $L_n=d_n2^{p_n-1}$, где ${d_n\in(0,1/2)}$. Рассмотрим на ${M_n=[0, 2^{-p_n})\times [0, 2^{p_n})}$ функцию  $f_n$, которая на прямоугольнике $[0, 2^{-p_n})\times [0, L_n)$ равна  ${a_n>0}$, на прямоугольнике ${[0, 2^{-p_n})\times [2^{p_n-1}, 2^{p_n-1}+L_n)}$ равна  $-a_n$, а в остальных точках $M_n$ равна нулю. Нетрудно проверить, что для любого ${x_n\in[0,2^{-p_n})}$
$$
\int_0^{2^{p_n}}f_n(x_n,y_n)\,dy_n=0,\ \ \int_0^{2^{p_n}}|f_n(x_n,y_n)|\,dy_n=2a_nL_n.
$$
Чтобы для функции ${f(x,y)=\sum\limits_{n=1}^\infty f_n(x_n,y_n)}$ не выполнялось условие~(3) из леммы~1, накладываем условие
$$
\sum_{n=1}^{\infty}a_nL_n=\sum_{n=1}^{\infty}a_nd_n2^{p_n}=\infty.
$$
При этом, чтобы ${f\in L_1(M)},$ достаточно сходимости ряда из $L_1(M_n)$-норм функций $f_n,$ т.е. ряда ${\sum\limits_{n=1}^{\infty}a_nd_n}.$ Мы же потребуем более сильное условие 
$$
\sum_{m=n+1}^\infty a_m=o(a_nd_n)\ \ \text{при}\ \ n\to+\infty.
$$ 

\begin{theorem}\label{Th2}
Пусть $f$ --- построенная выше функция  и $t_n=2^{p_n-2}.$ Тогда последовательность средних $A(f,t_n,(x,y))$ для потока над 2-одометром равномерно сходится к $0$ при $n\to\infty$ со скоростью $\mathcal{O}(a_nd_n)$.
\end{theorem}

{\bf Доказательство.}
Представляя 
$$
A(f, t_n, (x,y))=A\left(\sum_{m=1}^{n-1}f_m, t_n, (x,y)\right)+A(f_n, t_n, (x,y))+
A\left(\sum_{m=n+1}^{\infty}f_m, t_n, (x,y)\right),
$$
замечаем, что для всех ${(x,y)\in M}$ 
$$
A\left(\sum_{m=1}^{n-1}\, f_m, t_n, (x,y)\right)=0.
$$ 
Это следует из того, что промежуток интегрирования $[0,t_n)$ разбивается на конечное число промежутков длины $2^{p_m}$ для каждого ${m=1,...,n-1}.$ Это, действительно, будет так, поскольку $2^{p_n-2}$ делится на $2^{p_m}$. Отсюда возникает уточняющее условие на возрастающую последовательность $p_n:$
$$
p_{n+1}\geq p_n+2.
$$
А на интервале длины $2^{p_m}$ интеграл от $f_m$ по переменной $y_m$ равен нулю. 

Теперь обратим внимание на распределение значений функции 
$$
{A_n(x,y)=A\left(f_n, t_n, (x,y)\right)}.
$$
На множестве меры ${\frac 1 2 - d_n}$ функция $A_n(x,y)$  равна $0$, на  множестве меры ${1/ 4 - d_n/2}$ она равна $2a_nd_n$ (максимальное значение) и  на  множестве с такой же мерой $A_n(x,y)$ равна $-2a_nd_n$ (минимальное значение). 
Остальные промежуточные значения меняются  линейным образом.  

Покажем теперь, как оцениваются значения функции $A\left(\sum_{m=1}^\infty \, f_m,  t_n, (x,y)\right)$. Будем считать, что для всех $m>n$
выполняется условие 
$$
2^{p_n-2}<d_m2^{p_m-1}.
$$ 
Поскольку теперь длина промежутка интегрирования, равная $2^{p_n-2},$
существенно меньше $L_m$ при $m>n,$ то наибольшее по модулю значение среднего $A(f_m, t_n, (x,y))$ будет $a_m.$ Таким образом,
$$
\left|A\left(\sum_{m=n+1}^\infty\, f_m, t_n , (x,y)\right)\right|\leq\sum_{m=n+1}^\infty\,A(|f_m|, t_n , (x,y))\leq\sum_{m=n+1}^\infty a_m=o(a_nd_n).
$$
Собирая оценки вместе, получим для всех ${(x,y)\in M}$
$$
A(f, t_n, (x,y))=\mathcal{O}(a_nd_n)\ \ \text{при}\ \ n\to\infty.
$$
Теорема~\ref{Th2} доказана.

Значения $A\left(f, t_n ,(x,y)\right)$  на множестве меры близкой к $1/2$ по модулю близки к $2a_nd_n$. А для большинства остальных точек   значения функции  $A\left(f, t_n ,(x,y)\right)$ асимптотически являются $o(a_nd_n)$. Таким образом, нельзя получить оценку вида $o(a_nd_n).$  При этом для заданных последовательностей $a_n$ и $d_n$ мы можем выбирать  $t_n\to\infty$, растущую сколь угодно быстро.  
Таким образом, мы получаем сколь угодно медленную сходимость средних Биркгофа,  причем распределение значений этих средних таково, что на   почти половине пространства функция принимает значения асимптотически бесконечно большие  по сравнению со значениями, которые функция принимает на оставшейся части пространства. 

\section{Медленные сходимости средних Биркгофа в общем случае}

Мы показали, как реализуется медленная сходимость средних для для потока  над 2-одометром. Аналогичный эффект можно получить в общем случае.  Мы отчасти повторяем предыдущее построение,  при выборе   $t_{n+1}$ и $f_{n+1}$
применяем  специальное прямоугольное представление потока Рудольфа и эргодическую теорему Биркгофа.   Момент $t_{n+1}$ выбирается таким, что   средние  $A\left(f_1+\dots +f_n, t_{n+1} ,x\right)$ для большинства $x$ чрезвычайно малы по сравнению с $a_{n+1}$.

Для выбора $f_{n+1}$ используется теорема  о специальном прямоугольном представлении~(см., например,~\cite[глава~11,\S4]{KSF}). По  теореме Рудольфа для апериодического эргодического потока $T_t$ на вероятностном пространстве~${(X,m)}$ и любых положительных чисел $p,q$ таких, что $p/q$ иррационально, существует специальное представление потока с функцией, принимающей два значения: $p$  и $q.$ Таким образом,  фазовое пространство потока можно разбить на два измеримых множества, которые мы отождествляем с прямоугольниками. Первое множество состоит из отрезков траекторий длины $p$, а второе --- длины $q$.   

\begin{theorem}\label{Th4}
Пусть $T_t$ --- апериодический эргодический поток на вероятностном пространстве~${(X,m)}$ и ${\varphi(t)\to+0}.$ Тогда найдется функция ${f\in L_1(X,m)}$ с нулевым средним такая, что для п.в. ${x\in X}$
$$
\limsup_{t\to+\infty}\frac{1}{\varphi(t)}|A(f,t,x)|=+\infty.
$$
\end{theorem}

{\bf Доказательство.} Без ограничения общности можно считать, что $\varphi(t)$ монотонно стремится к нулю. Функцию ${f\in L_1(X,m)}$  ищем в виде 
$$
f(x)=\sum_{n=1}^\infty f_n(x_n,y_n),
$$
где $f_n$ строится следующим образом. В представлении Рудольфа   
выбираем параметры высот прямоугольников: $p=h_n$ и $q=h_n+\varepsilon_n$. При этом
$$
h_nc_n+(h_n+\varepsilon_n)d_n=1,
$$
и ${h_n\to+\infty},$ ${\varepsilon_n\to0}.$

Функция $ f_n(x_n,y_n)$, будет зависеть только от высоты $y_n$.  
А именно, пусть ${f_n=a_n}$ при ${0\leq y_n<h_n/4};$  ${f_n=-a_n}$ при ${h_n/2\leq y_n< 3h_n/4}$, а на остальной части $f_n$ равна 0 (см. Рис. 3).

\begin{figure}
	\begin{center}
		\includegraphics[scale=0.5]{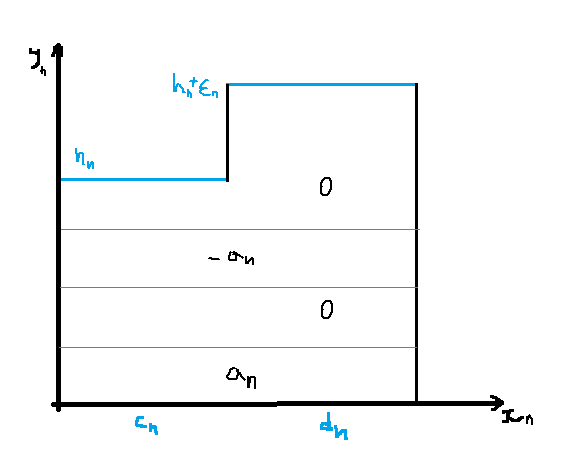} 
		
		{\caption{\small Представление Рудольфа с параметрами $h_n$ и $h_n+\varepsilon_n$}}
	\end{center}
\end{figure}

Всю эту «картинку» обратным изоморфизмом измеримо реализуем в исходном пространстве ${(X,m)}$ (точке $x$ соответствует пара $(x_n,y_n)$). Обозначение для функции оставляем таким же. Пусть ${\delta_n\to+0}$ и ${t_n=h_n\delta_n}.$ Для управления эргодическими средними зададим асимптотические условия на все последовательности: 
$$
\sum_{n=1}^\infty\varepsilon_nd_n+4\delta_n<1\eqno(I)
$$
$$
a:=\sum_{n=1}^\infty a_n<+\infty\ \ \text{и}\ \ \sum_{m=n+1}^\infty a_m=\alpha_na_n=o(a_n);\eqno(II)
$$  
$$
\frac{h_1+h_2+...+h_{n-1}}{h_na_n\delta_n}\to0\ \ \ \text{и}\ \ \ \frac{a_n\alpha_n}{\varphi(t_n)}\to+\infty.\eqno(III)
$$
По условию $(I)$ мы задаем последовательности ${\varepsilon_n}$ и ${\delta_n}.$ Далее, по условию~$(II)$ определяем $a_n$ и, самое важное для управления эргодическими средними, задаем рост $h_n$ условиями $(III).$ 

Как и ранее, представим  $A(f,t_n,x)$ в виде суммы трех  слагаемых:
$$
A(f, t_n, x)=A\left(\sum_{m=1}^{n-1}f_m, t_n, x\right)+A(f_n, t_n, x)+
A\left(\sum_{m=n+1}^{\infty}f_m, t_n, x\right).
$$
Для третьего слагаемого из условия~$(II)$ для всех ${x\in X}$ имеем
$$
\left|A\left(\sum_{m=n+1}^{\infty}f_m, t_n, x\right)\right|\leq\sum_{m=n+1}^\infty a_m=o(a_n).
$$
Для первого слагаемого из первой части условия~$(III)$, в виду зануления интеграла на целых отрезках траекторий, для всех ${x\in X}$ имеем
$$
A\left(\sum_{m=1}^{n-1}|f_m|, t_n, x\right)\leq\frac{2}{t_n}\sum_{m=1}^{n-1}a_mh_m\leq
2a\frac{h_1+h_2+...+h_{n-1}}{h_n\delta_n}=o(a_n).
$$
Для второго слагаемого нетрудно посчитать, что оно принимает значение $a_n$ на множестве меры $(c_n+d_n)h_n(1/4-\delta_n)$, на другом множестве такой же меры --- значение $-a_n,$ и на множестве меры не меньше $2(c_n+d_n)h_n(1/4-\delta_n)$ принимает нулевое значение. Остальные значения по модулю не превосходят $a_n.$ 

Таким образом, $A(f,t_n, x)=\mathcal{O}(a_n)$ для всех $x\in X.$ При этом, на некотором множестве $\mathcal D_n$ меры 
$$
m(\mathcal{D}_n)=2(c_n+d_n)h_n(1/4-\delta_n)=2(1/4-\delta_n)(1-\varepsilon_nd_n)>1/2-2\delta_n-\frac{\varepsilon_nd_n}{2}
$$ 
$$
|A(f,t_n,x)|=a_n+o(a_n);
$$
а на некотором множестве $\mathcal{E}_n,$ с такой же оценкой на меру
$$
m(\mathcal{E}_n)>1/2-2\delta_n-\frac{\varepsilon_nd_n}{2},
$$
будет
$$
|A(f,t_n,x)|=o(a_n).
$$
Пусть ${\mathcal{F}=\bigcap_{n=1}^\infty(\mathcal{E}_n\cup\mathcal{D}_n)},$ тогда, учитывая~$(I),$ ${m(\mathcal{F})>1-\sum_{n=1}^\infty\varepsilon_nd_n+4\delta_n>0}.$

Для всех $x\in\mathcal{F}$ для рассматриваемых средних имеем 
$$
\text{либо}\ \  |A(f,t_n,x)|=a_n(1+\alpha_n), \ \ \text{либо}\ \ |A(f,t_n,x)|=\alpha_na_n.
$$
В любом случае, для таких $x$ из второго условия в~$(III)$ получим
$$
\frac{ |A(f,t_n,x)|}{\varphi(t_n)}\to+\infty,
$$
т.е. для $x\in\mathcal{F}$
$$
\limsup_{t\to+\infty}\frac{1}{\varphi(t)}|A(f,t,x)|=+\infty.
$$
Для п.в. $x\in X\setminus \mathcal{F}$ будет выполняться такое же соотношение ввиду закона нуля или единицы для скорости сходимости~\cite{KaPoS2022}.

\section {О медленной сходимости временных средних для обмотки тора и непрерывной функции}

Оказывается, что усреднение непрерывной  функции вдоль траекторий  гладкого потока  можно совместить с  эффектом  медленной сходимости средних.  Для этого  мы подходящим образом будем сглаживать  функции $f_n$, рассмотренные  выше.  Пусть $T_t(x,y)=(\{x+t\},\{y+ct\}), {(x,y)\in M}$ --- эргодическая обмотка тора вдоль 
вектора $(1,c)$, где $c$ --- иррациональное число, которое приближается рациональными дробями  $p_n/q_n$ так, что 
$$
0<c-\frac{p_n}{q_n}=\beta_n\frac{1}{q_n^2}=o\left(\frac{1}{q_n^2}\right).
$$

Рассматривая обмотку тора вдоль вектора ${(1,p_n/q_n)},$ отождествим  тор~$M$  с узким параллелограммом~$M_n$ со стороной длины $\sqrt{p_n^2+q_n^2}$, параллельной этому вектору,  и стороной длины $1/q_n$, параллельной оси OX.  Внутри этого параллелограмма, рассмотрим параллелограмм~$R_n,$ направленный вдоль потока $T_t.$
\begin{figure}
	\begin{center}
		\includegraphics[scale=0.5]{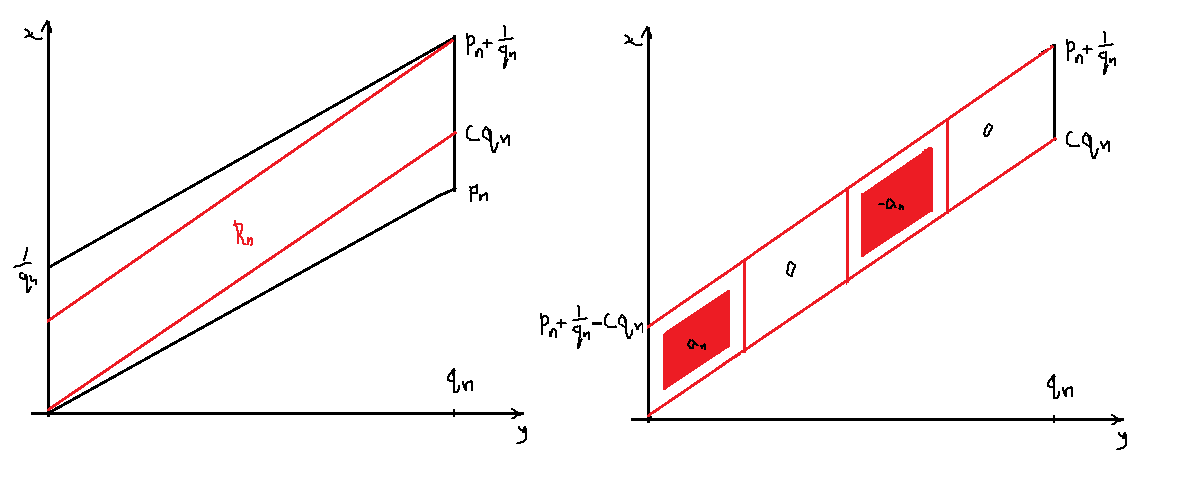} 
		
		{\caption{\small Параллелограмм $R_n$ и определение функции $f_n$ на нём}}
	\end{center}
\end{figure}
Нетрудно видеть, что площадь параллелограмма~$R_n$ равна
$$
1-q_n^2\left(c-\frac{p_n}{q_n}\right)=1-\beta_n.
$$

На $M_n$ выберем координаты $x_n\in[0,1/q_n]$ вдоль оси OX и $y_n\in[0,q_n\sqrt{1+c^2}]$ вдоль вектора $(1,c).$
Построим гладкую функцию $f_n(x_n,y_n)$, которая по модулю не превосходит ${a_n>0}$, с условием ${\sum_n a_n<\infty}$. А именно, разделим $R_n$ на четыре равные по площади части. На первой четверти $f_n$ равна $a_n$ на множестве, отделенным от границы и площади близкой к $1/4.$ На оставшейся части этой четверти $f_n$ гладко убывает до нуля на границе. На втором и четвертом параллелограмме $f_n$ равна 0. На третьей четверти значения антисимметричны значениям на первой четверти (см. Рис.~3).

Вне прямоугольника $R_n$ и на его границе пусть $f_n$ равна $0$. Также считаем, что 
$$
\int_0^{q_n\sqrt{1+c^2}}f_n(x_n,y_n)\,dy_n=0
$$ 
для каждого ${x_n\in[0,1/q_n]}.$ 

В результате получим непрерывную функцию $f(x,y)=\sum\limits_{n=1}^\infty f_n(x_n,y_n)$. 

Пусть ${\varphi(t)\to+0}$ при $t\to+\infty.$ Положим $t_n=\delta_nq_n\sqrt{1+c^2},$ а также пусть выполняются условия $(II)$ и $(III),$ где вместо $h_n$ берем $q_n.$ 

\begin{theorem}\label{Th3}
	Пусть $f$ --- построенная выше функция. Тогда эргодические средние $A(f,t,(x,y))$ для обмотки тора вдоль вектора ${(1,c)}$ удовлетворяют для п.в. ${(x,y)\in M}$ соотношению
	$$
	\limsup_{t\to+\infty}\frac{1}{\varphi(t)}|A(f,t,(x,y))|=+\infty.
	$$
\end{theorem}

{\bf Доказательство.} Временн\'{о}е среднее представим в виде  
$$
A(f,t_n,(x,y))=A\left(\sum_{m=1}^{n-1}f_m, t_n, (x,y)\right)+A(f_n,t_n,(x,y))+A\left(\sum_{m=n+1}^\infty f_m, t_n, (x,y)\right).
$$ 
Для всех ${(x,y)\in M}$ имеем очевидную оценку для третьего слагаемого 
$$
\left|A\left(\sum_{m=n+1}^\infty f_m, t_n, (x,y)\right)\right|\leq \sum_{m=n+1}^nA(|f_m|, t_n, (x,y))\leq\sum_{m=n+1}^\infty a_m=o(a_n).
$$
Чтобы оценить  первое слагаемое в указанной выше сумме, воспользуемся занулением интегралов вдоль траектории потока. На  большей части отрезков траекторий, вдоль которых берется интеграл, вклад в интеграл нулевой в силу выбора функций $f_m$.
Нам нужно  оценить  только значения интегралов  вдоль начальной и конечной частей отрезков траектории. А именно, для каждого $m=1,...,n$ и всех $(x,y)\in M$
\begin{multline*}
t_nA(f_m, t_n, (x,y))=\int_0^{q_m\sqrt{1+c^2}-y_m}f_m(T_s(x_m,y_m))\,ds+
\underbrace{\int_0^{q_m\sqrt{1+c^2}}f_m(x^{(1)}_m,s)\,ds}_{=0}+...\\
...+\underbrace{\int_0^{q_m\sqrt{1+c^2}}f_m(x^{(k-1)}_m,s)\,ds}_{=0}+\int_0^{y_m}f_m(x^{(k)}_m,s)\,ds.
\end{multline*}
Отсюда уже получаем   
\begin{multline*}
\left|A\left(\sum_{m=1}^{n-1}f_m, t_n, (x,y)\right)\right|\leq\sum_{m=1}^{n-1}\frac{2a_mq_m\sqrt{1+c^2}}{t_n}=\mathcal{O}\left(\frac{q_1+\dots+q_{n-1}}{\delta_nq_n}\right)=o(a_n).
\end{multline*}

Второе слагаемое ведет себя так же, как и в в теореме~\ref{Th4} для общего случая. А именно, на множестве, близком по мере к $1/2$, средние ${|A(f,t_n,(x,y))|=a_n+o(a_n)},$ а на другом множестве, также близком по мере к $1/2,$ средние ведут себя как  ${o(a_n)}.$  
Теорема~\ref{Th3} доказана.

\section{Заключительные замечания}

С тематикой  скоростей сходимости эргодических средних сопряжено большое количество нерешенных задач.

{\bf Сходимость средних и гладкость  усредняемой функции.} Построенная в теореме~\ref{Th3} функция является лишь непрерывной. Интересен вопрос о существовании гладких функций со сколь угодно медленной скоростью сходимости средних. Как показал Ковада~\cite{Kow} (см. также~\cite{Web}), начиная с некоторого показателя гладкости, зависящего от скорости аппроксимации иррационального числа $c$, скорость сходимости эргодических средних для обмотки тора будет максимальной $\mathcal{O}(1/t).$

\vspace{3mm}

{\bf Возможные распределения средних Биркгофа.}
Мы показали, как выбор подходящей усредняемой функции реализует диапазон 
скоростей сходимости от максимальной до сколь угодно медленной.
Эффекты замедления скоростей сходимости средних можно обнаруживать для широкого 
класса групповых действий~(см.~\cite{R}).  
Интерес представлет не только оценка скоростей сходимости, но и более общая задача  
о возможных распределениях значений средних Биркгофа. Благодаря выбору усредняемой функции  этими  распределениями  можно управлять. 
Например, пусть   мы хотим, чтобы распределение значений функции 
$A(f, t_n, x)$ для эргодического потока было сколь угодно близко 
к распределению значений, например, функции  $c_n s^2 , s\in[0, 1],$ 
относительно меры Лебега на $[0, 1].$ Если  числа $c_n$  
достаточно быстро стремятся к $+0,$  задача  решается методом, 
похожим на изложенный выше,  при подходящем  выборе  функций $f_m$.

\vspace{3mm}
{\bf Поиск оптимальных  весовых  распределений.}  Как уже отмечалось, при рассмотрении весовых усреднений возможно увеличение  диапазона быстрых скоростей сходимости средних~(см.~\cite{Izv}, а также ссылки в этой работе). В связи с этим возникает ряд новых задач. Сформулируем следующие частные случаи.

Пусть  $f$ ---  функция с  нулевым средним заданной гладкости $k\geq1$. 
Для эргодического сдвига $T$ на торе рассмотрим всевозможные выпуклые суммы вида
$$
A_w(f,N,x)=\sum _{n=1}^Nw_{n,N}f(T^nx),\ \ 
$$
где $w=(w_{1,N},\dots,w_{N,N})$ --- вероятностный вектор.  

При каких  вероятностных распределениях коэффициентов $w_{n,N}$ при заданном достаточно большом $N$  нормы  $\|A_w\|$ (в пространствах $L_1,$ или  $L_2,$ или $L_\infty$) будут минимальны?

Аналогичный вопрос ставится для непрерывного времени и обмоток тора.
Для средних
$$
A_w(f,t,x) =\int_0^t w(s,t) f(T_s x)\,ds
$$
где $w\geq 0$  и   $\int_0^t w(s,t)ds =1$ для каждого $t>0,$
ищем $\inf_{w}\|A_w\|.$         

Даже в случае вращения окружности такие задачи представляются нетривиальными.

\text{}

Работа Подвигина И.В. выполнена в рамках государственного задания ИМ СО РАН
(проект № FWNF-2026-0022).

\end{document}